\documentclass{amsart}

\usepackage{graphicx}

\graphicspath{ {./} }

\usepackage{amssymb,amsfonts}
\usepackage[all,arc]{xy}
\usepackage{enumerate}
\usepackage{mathrsfs}
\usepackage{enumitem}
\usepackage{amssymb}
\usepackage[stable]{footmisc}

\newcommand{\bbN}{\mathbb{N}}

\newcommand{\bbR}{\mathbb{R}}
\newcommand{\bbZ}{\mathbb{Z}}

\newcommand\Z[1]{\mathbb{Z}/#1\mathbb{Z}}

\newtheorem{thm}{Theorem}[section]
\newtheorem{cor}[thm]{Corollary}
\newtheorem{prop}[thm]{Proposition}

\newtheorem{rmrk}[thm]{Remark}

\theoremstyle{definition}
\newtheorem{defn}[thm]{Definition}

\theoremstyle{remark}
\newtheorem{rem}[thm]{Remark}

\makeatletter
\let\c@equation\c@thm
\makeatother
\numberwithin{equation}{section}

\begin{document}

\title{A bottom-up approach to Hilbert's Basis Theorem}

\author{ Marc Maliar }

\date{\today}

\begin{abstract}
In this expositional paper, we discuss commutative algebra---a study inspired by the properties of integers, rational numbers, and real numbers. In particular, we investigate rings and ideals, and their various properties. After, we introduce the polynomial ring and the fundamental relationship between polynomials and sets of points. We prove some results in algebraic geometry, notably Hilbert's Basis Theorem. 
\end{abstract}

\maketitle

\tableofcontents

\section{Introduction}  
Mathematics often begins with concrete problems, for example finding the length of the hypotenuse of a triangle. But as a mathematician investigates such problems, they begin to abstract away various amounts of details, uncovering patterns and defining properties.

In this expositional paper, we take this approach to commutative algebra. 

In Section $2$, we attempt to find the patterns and defining properties hidden in algebra. In particular, we investigate sets of numbers such as the integers, rational numbers and real numbers. We define a ring, the generalization of the properties of such sets of numbers. We realize that there exist special sets in the integers called ideals---named in this way because they behave in an ``ideal" way. Finally, we explore properties of these rings and ideals.

But we realize that polynomials behave as integers do: we can add and multiply polynomials, identify ``prime" polynomials, and find ideals. So we begin Section $3$ by defining the polynomial ring. Unlike the integers, though, polynomials have an interesting relationship with sets of points in a field. And so the rest of section $3$ is devoted to applying commutative algebra to this relationship between sets of points and polynomials. For example, we establish that, given a set of polynomials $S$ in two variables (or any number of variables actually), the set of points $X=(x,y)$ such that $f(x,y)=0$ for all $f \in S$ can be written as the intersection of finitely many sets of points $X_1,\dots,X_n$, where $X_i=(x,y)$ such that $f_i(x,y)=0$ for some polynomial $f_i$. This result is a corollary of Hilbert's Basis Theorem (Theorem \ref{noetherian}). 

By the end, we are able to use abstractions to prove nontrivial theorems about sets of points and polynomials. 
\section{Rings and ideals}
We begin this section by exploring rings. What is a ring? Consider the set of integers, $\mathbb{Z}$. Recall that this fundamental set comes equipped with two binary operations (the usual addition and multiplication). We generalize the properties of this set and arrive at the definition of a ring.

We then consider special subsets of a ring, called ideals. We motivate this study by again considering $\mathbb{Z}$. In particular, we look at sets of integers that are divisible by a certain integer, and we note that they possess interesting properties. We then define ideals by these interesting properties. In the rest of the section, we discover other interesting properties of ideals and investigate related definitions.
\subsection{Definitions of a ring}
\begin{defn}\label{ring}
A \emph{ring} $R$ is a set with two binary operations addition and multiplication, denoted by $+$ and $*$ respectively. These operations have the following properties:
\begin{enumerate}
    \item Addition is associative and commutative; there exists an additive identity (denoted by $0$); each element in the ring has a (unique) additive inverse.
    \item Multiplication is associative; there exists a multiplicative identity (denoted by $1$).
    \item Multiplication is distributive over addition. That is, $r_1*(r_2+r_3)=r_1*r_2+r_1*r_3$ and $(r_1+r_2)*r_3=r_1*r_3+r_2*r_3$ hold for $r_1,r_2,r_3 \in R$. 
\end{enumerate}
\end{defn}
Why such a definition? The ring definition comes from trying to generalize numbers (integers, rational numbers, and real numbers) and their various properties. Notice that the aforementioned sets of numbers with their usual operations satisfy these properties.

We consider commutative rings: rings in which multiplication is commutative. From now on, we consider only commutative rings; any ring mentioned is commutative. Notice that the integers, rational numbers, and real numbers are commutative rings. 

While each element in a ring has an additive inverse, it does not necessarily have a multiplicative inverse. In addition, $1$ can be equal to $0$. In this case, it is easy to see that the ring has only one element, $0$ (and that this set of one element is a ring). 

Recall the classical definition of a field: a set with addition and multiplication satisfying the usual properties. Fields turn out to be valuable objects in algebraic geometry, so I mention them here. Also, in the next section, we prove some results about commutative rings that are fields, so I introduce the following (equivalent) definition:
\begin{defn}
A \emph{field} is a commutative ring $k$ where $1 \neq 0$ and every element has a multiplicative inverse. 
\end{defn}
\subsection{Definition of an ideal}
Consider the ring of integers; now, consider the set of all elements divisible by $3$ (a certain subset of this ring). This set has some interesting properties.

First, if you add any two elements divisible by $3$, you end up with an integer divisible by $3$. Second, if you multiply an element in the set with any other integer, you end up with an integer divisible by $3$. Ideals are defined by these two properties.
\begin{defn}\label{ideal}
An \emph{ideal} $I$ is a subset of a ring $R$ that, inheriting the addition and multiplication operations from $R$, satisfies the following properties:
\begin{enumerate}
    \item Closure with respect to addition (if $i_1,i_2 \in I$ then $i_1+i_2\in I$); associativity with respect to addition; additive identity is in the ideal ($0 \in I$); each element's additive inverse is also in the ideal (if $i \in I$, $-i \in I$).
    \item If $r \in R$ and $i \in I$, $r*i \in I$.
    \item If $r \in R$ and $i \in I$, $i*r \in I$. 
\end{enumerate}
\end{defn}
\begin{rmrk}
Since we work with commutative rings, it is only necessary to prove either $(2)$ or $(3)$. In non-commutative rings, an ideal that satisfies $(2)$ but not necessarily $(3)$ is a left-ideal, and an ideal that satisfies $(3)$ but not necessarily $(2)$ is a right-ideal.
\end{rmrk} 
Note that we do not need to explicitly prove all of these statements. Indeed, proving closure with respect to addition, existence of additive inverse, and $(2)$ (or $(3)$ as discussed) is sufficient. 
\begin{enumerate}[label={(\alph*)}]
    \item Associativity is inherited from the ring definition.
    \item For $i \in I$, $i*(-1)=-i \in I$. 
\end{enumerate}
An ideal cannot be empty because it must contain $0$. 

Let us find other ideals of the rings we have mentioned (integers, rational numbers, and real numbers). To begin, note that, for any ring, both $\{0\}$ and $R$ are ideals. 

But are these the only ideals? Consider first the integers.
\begin{prop}
For all $z \in \bbZ$, define $I_z=\{z*r \mid r \in \bbZ\}$. $I_z$ is an ideal. These are the only ideals of $\bbZ$.
\begin{proof}
The proof that $I_z$ is an ideal follows from my previous remark. So then suppose there exists an ideal $I$ that cannot be written as an $I_z$. $I \neq \{0\}$ because then $I=I_0$, so there exists at least one nonzero integer in $I$. Consider the set of nonzero positive integers in $I$. This set must be nonempty; an empty set would imply that the nonzero integer $z$ we found earlier is negative, but $z*(-1) \in I$---a contradiction. By the Well-Ordering Principle, find the smallest positive nonzero integer $z \in I$. 

Because $I \neq I_z$, there must exist an element $z' \in I$ that is not in the form $z*r$ for any $r \in \bbZ$. Suppose it is positive (if it is negative multiply it by $-1$ and it will still be in the ideal). We know that $z'>z$ because $z$ is the smallest positive integer. Now, observe that $z'+(-z)$ is still in the ideal and also will not be in the form $z*r$. Repeat until the result is less than $z$. The result will be in $I$, positive and nonzero (because it is not in the form $z*r$). Thus, we have found a smaller positive nonzero integer in $I$ that is smaller than $z$---a contradiction. 
\end{proof}
\end{prop}
In other words, the only ideals of the ring of integers are multiples of an integer. Now, consider rational numbers and real numbers. It turns out that we can make a statement about fields in general:
\begin{prop}
The only ideals of a field $k$ are the field itself and the zero ideal. The converse is also true: if a commutative ring has exactly two ideals, the zero ideal and the entire ring, then the commutative ring is a field. 
\begin{proof}
First, we prove that that the only ideals of $k$ are the field itself and the zero ideal. Let $I$ be an arbitrary ideal.
\begin{enumerate}[label={Case \arabic*:}, wide=15pt, leftmargin=*]
    \item $I$ contains no non-zero elements. Then, $I=\{0\}$.
    \item $I$ contains a non-zero element. Let $i$ be this element. Then, because $k$ is a field, there exists a multiplicative inverse of $i$, denoted by $i^{-1}$. By the property of ideals, $i*i^{-1}=1 \in I$. Any element $j \in k$ is in the ideal now because $1*j=j \in I$, so then $I=k$. 
\end{enumerate}
$\indent$ Second, we prove that if the only ideals of a ring $R$ are the ring and the zero ideal, the ring is a field. Let $r$ be a non-zero element in our commutative ring $R$. Construct a set $I=\{rq \mid q \in R\}$. It is easy to check that this is an ideal. Because this ideal contains a nonzero element $r$, it cannot be the zero ideal, so by our given, the ideal must be the entire ring. Then, there exists an element $r'$ such that $r(r')=1$ since $1$ is in our ideal. Thus, $r$ has an inverse. 
\end{proof}
\end{prop}
From this proof, we can also see the following statement:
\begin{cor}\label{ringideal}
Let $I$ be an ideal of $R$. Then, $1 \in I \iff I=R$.
\end{cor}
\subsection{Another definition of an ideal}
Consider the following definition:
\begin{defn}
A \emph{ring homomorphism} is a map $\phi$ from one ring $R$ to another ring $S$ with the following properties:
\begin{enumerate}
    \item $\phi(r_1+r_2)=\phi(r_1)+\phi(r_2)$,
    \item $\phi(r_1*r_2)=\phi(r_1)*\phi(r_2)$, and
    \item $\phi(1_R)=1_S$, where $1_R, 1_S$ are the multiplicative identities of $R$ and $S$ respectively.
\end{enumerate}
\end{defn}
We introduce an equivalent definition of the ideal using ring homomorphisms:
\begin{defn}\label{homomorphism}
An \emph{ideal} $I$ is a subset of a ring $R$ that is the kernel of a ring homomorphism. That is, a set $I\subset R$ is an ideal if there exists a ring homomorphism $\phi\colon R\to S$, where $S$ is some other ring, such that $\phi^{-1}(0)=I$. 
\end{defn}
\begin{prop}
Definition $\ref{ideal}$ is equivalent to Definition $\ref{homomorphism}$.
\begin{proof}
First, we show Definition $\ref{homomorphism}$ from Definition $\ref{ideal}$. We need to find a ring $S$ and a function $\phi$ as described. We need $\phi$ such that for all $i \in I$, $\phi(i)=0$. Furthermore, we need $\phi$ such that adding $i \in I$ to any $r \in R$ will not change the computed result (i.e., $\phi(r)=\phi(r+i)$). 

Let $r$ be an element in $R$. Define $\phi\colon r\mapsto X_r:=\{r+i\ :\ i\in I\}$. Define $R/I=\{X_r \mid r \in R\}$. 

We now define an equivalence relation between elements in $R$. We say that $r \sim q$ for $r,q \in R$ if and only if $\phi(r)=\phi(q)$. It follows from properties of an ideal that $\sim$ is an equivalence relation. Notice that, under this definition, $\phi(i)=\phi(0) \; \forall \; i \in I$. 

Define addition and multiplication as follows: for $\phi(r),\phi(q) \in R/I$, let $\phi(r)+\phi(q)=\phi(r+q)$ and $\phi(r)*\phi(q)=\phi(r*q)$. Both of these are well-defined. Consider, for example, addition. Let $r_1,r_2,q_1,q_2$ be elements in $R$ with $r_1 \sim r_2$ and $q_1 \sim q_2$. Let $i,j$ be elements in $I$ such that $r_1=r_2+i$ and $q_1=q_2+j$. Then,

$$\phi(r_1)+\phi(q_1)=\phi(r_1+q_1)=\phi(r_2+q_2+i+j)=\phi(r_2+q_2)=\phi(r_2)+\phi(q_2).$$

We conclude that equivalent inputs implies equivalent outputs under our addition operation. Multiplication follows similarly as well.

Under these operations, $R/I$ forms a ring. To finish the proof, define $S=R/I$ and $\phi$ as defined above. By our definition of $\phi$, we already know that $\phi(r_1+r_2)=\phi(r_1)+\phi(r_2)$ and $\phi(r_1*r_2)=\phi(r_1)*\phi(r_2)$. It suffices to show that $\phi(1_R)=1_S$ where $1_S$ is the multiplicative identity of $R/I$. The proof follows from the definition.

Second, we show Definition $\ref{ideal}$ from Definition $\ref{homomorphism}$. As discussed, we must prove only three statements:
\begin{enumerate}[label={(\alph*)}]
    \item Closure with respect to addition: let $i_1,i_2 \in I$. We know $\phi(i_1+i_2)=\phi(i_1)+\phi(i_2)=0+0=0$. Then, $i_1+i_2 \in I$ because $\phi(i_1+i_2)=0$. 
    \item Existence of additive inverse: $\phi(0)=0$ so $0 \in I$. 
    \item If $i \in I$ and $r \in R$ then $i*r \in I$: let $i \in I, r \in R$. $\phi(i*r)=\phi(i)*\phi(r)=0*\phi(r)=0$. Then, $i*r \in I$ because $\phi(i*r)=0$. 
\end{enumerate}
\end{proof}
\end{prop}
Consider the ideals of $\bbZ$ we have found. What is the ring $S$ and ring homomorphism $\phi$ that we should be able to find, according to our new definition of an ideal?
\begin{enumerate}
    \item $I_z=\{nz \mid z \in \bbZ\}$: we use $S=\Z{n}$ ($S$ is the set of integers modulo $n$---a ring) and $\phi$ where $\phi(r)=r \mod n$. The set of elements that map to $0$ is $\{nz \mid z \in \bbZ\}$.
    \item $I=\{0\}$: we use $S=R$ and $\phi$ where $\phi(r)=r$. The only element that maps to $0$ is $0$.
    \item $I=\bbZ=R$: we use $S=\{0\}$ and $\phi$ where $\phi(r)=0$. All elements map to $0$. 
\end{enumerate}
\subsection{What is in an ideal?}
Consider an ideal $I$ of ring $R$. How can we imagine what is in this ideal?

Let us try to come up with an expression that can represent any element in the ideal. Let $I$ be an ideal and let $\{i_\mu\}$ be a subset of $I$ indexed by $\mu$. Consider the expression $\sum_\mu i_\mu*r_\mu$ where $i_\mu$ is an element of the ideal and $r_\mu$ is an element in $R$. Can we write any element in $I$ in this form for some $r_\mu$'s? If so, we denote the ideal $I$ by
$$
I=(\{i_\mu\})
$$
and say that it is generated by $\{i_\mu\}$. For example, $I=R=(1)$ (as we saw in Corollary \ref{ringideal}). So, now we have come up with a way to imagine elements in our ideal.

Some things to emphasize: in rings that extend outwards infinitely (like the integers and the ring of polynomials introduced later), non-zero ideals grow large up to arbitrary finite size. For example, with the ideal of integers divisible by $3$, we can keep multiplying and finding larger and larger elements. 

In addition, note that just because we can write an element $i \in I$ as a sum of $i_\mu*r_\mu$'s doesn't mean that for any two elements $a,b \in R$ such that $a*b=i$, either $a$ or $b$ in $I$. An ideal in which this is true for every $i \in I$ is called a prime ideal:
\begin{defn}\label{primeideal}
A \emph{prime ideal} $I$ of ring $R$ is an ideal that satisfies the following properties:
\begin{enumerate}
    \item If $k \in I$ and $i*j=k$, then either $i \in I$ or $j \in I$. 
    \item $I$ is not equal to the ring $R$.
\end{enumerate}
\end{defn}
Consider the ideal of integers divisible by $3$ in the ring of integers: this ideal is a prime ideal because we can write every element as $3$ multiplied by an integer, and $3$ is in our ideal. However, the ideal of integers divisible by $6$ does not have this property; for example $6=2*3$ but neither $2$ nor $3$ are divisible by $6$. The generalization of these statements corresponds to the following statement:
\begin{prop}
In $\bbZ$, ideals generated by a prime number are prime. The zero ideal is also prime. All other ideals are not prime.
\begin{proof}
Let $I$ be an ideal generated by a prime natural number $p$. Suppose $I$ is not prime. Then, there exists an element $k \in I$ and $i,j \in I$ such that $i*j=k$ but $i \not\in I$ and $j \not\in I$. This means that $k$ is divisible by $p$ but neither $i$ nor $j$ are divisible by $p$, which is impossible.

The zero ideal is prime because no integer divides $0$. 

We already showed that all ideals of $\bbZ$ are generated by one element, and so it is only necessary to look at ideals generated by nonprime elements. Let $I$ be an ideal generated by a nonprime natural number $n$. Because $n$ is not prime, $n=\ell*m$ for some $\ell,m$ not divisible by $n$, so $\ell$ and $m$ are not in $I$; then, $I$ is not prime. The ideal $I$ generated by $-n$ is the same as the one generated by $n$ so the same logic applies. 
\end{proof}
\end{prop}
We have an alternate definition of prime ideal as well. This definition relates to Definition $\ref{primeideal}$ similarly to how Definition $\ref{homomorphism}$ relates to Definition $\ref{ideal}$.
\begin{defn}\label{primeidealhomomorphism}
A \emph{prime ideal} $I$ of a ring $R$ is the kernel of a ring homomorphism $\phi \colon R \rightarrow S$, where $S$ is a ring that is not equal to $\{X_0\}$ and has no zero-divisors (no elements $s$ such that for some $t \in S$, $s*t=0$). 
\end{defn}
\begin{prop}
Definition $\ref{primeideal}$ is equivalent to Definition $\ref{primeidealhomomorphism}$.
\begin{proof}
First, we prove that Definition $\ref{primeideal}$ is equivalent to Definition $\ref{primeidealhomomorphism}$. 

Let $I$ be a prime ideal by Definition $\ref{primeideal}$. Suppose that Definition $\ref{primeidealhomomorphism}$ does not hold. By Definition $\ref{homomorphism}$, there must exist a set $S$ and $\phi$ such that $I$ is the kernel of $\phi$. Then, take $S$ and $\phi$ that Definition $\ref{primeidealhomomorphism}$ refers to. Consider two cases:

\begin{enumerate}[label={Case \arabic*:}, wide=15pt, leftmargin=*]
    \item $S=\{X_0\}$. $I$ is not equal to $R$. Then, there exists an element $r$ in $R$ that is not in $I$, and $\phi(r) \neq X_0$. 
    
    \item $S$ has a zero-divisor. Let $s$ be a zero-divisor that is not equal to $0$ with a corresponding $t$ that is not equal to $0$ such that $s*t$ is equal to $0$. By construction, $\phi$ is surjective, so there exist $r,q \in R$ such that $\phi(r)=s$ and $\phi(q)=t$. Then, 

    $$0=s*t=\phi(r)*\phi(q)=\phi(r*q),$$

    so $r*q \in I$. However, $r \not\in I$---otherwise $\phi(r)=0$ (and similar argument for $s$). Then, $I$ is not prime according to Definition $\ref{primeideal}$, so we have arrived at a contradiction.
\end{enumerate}

Second, we prove that Definition $\ref{primeidealhomomorphism}$ is equivalent to Definition $\ref{primeideal}$.  

For sake of contradiction, let $I$ be an ideal that is not prime by Definition \ref{primeideal}. By Definition $\ref{homomorphism}$, there exists a set $S$ and $\phi$ such that $I$ is the kernel of $\phi$. Consider two cases:
\begin{enumerate}[label={Case \arabic*:}, wide=15pt, leftmargin=*]
    \item $I=R$. Then, every element $i \in I$ maps to $0$, so $S=\{X_0\}$. 
    \item There exist elements $r,q \in R$ such that $r*q$ is in $I$ but neither $r$ nor $q$ is in $I$. Then,

    $$0=\phi(r*q)=\phi(r)*\phi(q),$$

    but neither $\phi(r)$ nor $\phi(q)$ are equal to $0$ because neither $r$ nor $q$ are in $I$. Then, both $\phi(r)$ and $\phi(q)$ are zero-divisors of $S$ so we have arrived at a contradiction.
 
\end{enumerate} 
\end{proof}
\end{prop}
There are many other interesting ways to characterize rings and ideals. See \cite[Chapter~1]{atiyah}.
\subsection{Finitely generated ideals}
Let us go back to our equation for elements in an ideal. Can we make some claim about the sets $i_\mu$ that generate $I$? 

Note that such a set must exist. If we choose $\{i_\mu\}=I$ then $I$ is generated by $\{i_\mu\}$. Let $i_\lambda$ be an element in $I$. Set $r_\mu=1$ for $\mu=\lambda$ and $r_\mu=0$ for $\mu \neq \lambda$.

But can we find other sets? Specifically, can we find a finite set $\{i_\mu\}$ that generates $I$? Such an ideal is said to be finitely generated. Furthermore, for a ring $R$, is every ideal $I \subset R$ finitely generated? Such a ring is said to be a \emph{Noetherian} ring. 

Consider the ring of integers. Since every ideal of $\bbZ$ consists of all multiples of some element, every ideal $I$ is finitely generated and $\bbZ$ is Noetherian. 

Later, I will give an example of a ring that is not Noetherian. 

There are two other alternate definitions of Noetherian rings: 
\begin{prop}
The following are equivalent:
\begin{enumerate}[label=\roman*]
    \item For every ideal $I \subset R$, there exist $x_1,\dots, x_n \in I$ such that $I=(x_1,\dots,x_n)$.
    \item Every chain of ideals $I_1 \subset \dots \subset I_k \dots$ eventually terminates such that $I_\ell=I_{\ell+1}=\dots$. 
    \item Every nonempty set of ideals of $R$ has an ideal that is not contained in any other. 
\end{enumerate}
\begin{proof}
(i) $\Rightarrow$ (ii): Take a chain of ideals, and let $I$ be the union of them all. For $I$, there exist $x_1,\dots,x_n$ that generate $I$. All of these $x_i$ should appear at some point in the chain of ideals. When they appear, they stay in the chain because each ideal contains all previous ideals. Because there are finitely many $x_i$, we can find an $I_\ell$ for some $\ell$ that finally contains all of them. It is clear that $I_{\ell}=I$ because they are generated by the same elements. To prove the chain terminates, let $k$ be a natural number bigger than $l$. $I_l$ is a subset of $I$ and therefore is a subset of $I_k$. $I_k$ is a subset of $I_l$ because we have a chain of ideals. Then, $I_l=I_k$, so the chain terminates.  \\
\indent \indent \; (ii) $\Rightarrow$ (iii): At the ideal at which the chain terminates, that ideal is not contained in any other. \\
\indent \indent \; (iii) $\Rightarrow$ (i): Let $I$ be an ideal and $\Sigma=\{J \subset I \mid J \text{ is a finitely generated ideal}\}$. By our given, $J$ has an element $J_0$ that is not contained in any other. Suppose that $J_0 \neq I$. Then, there exists $x \in I, x \not\in J_0$. Let $J_0'$ be the ideal generated by all the finitely many generators in $J_0$ and also $x$. $J_0'$ is still finitely generated and contains $J_0$, contradicting our assumption that $J_0$ is not contained in any other. Then, $I=J_0$ is finitely generated. 
\end{proof}
\end{prop}
\section{Algebraic varieties and the polynomial ring}
We begin our study of a polynomial ring, points, and the relationship between them. We motivate our study first by introducing the plane curve, a mathematical object studied heavily in the $19$th century. Next, we formally define polynomials and the fundamental functions ``taking the variety of" and ``taking the ideal of" that map polynomials to points and vice-versa, respectively. After, we investigate how these two functions serve as inverses to each other. We then attempt to understand the nature of plane curves (and algebraic sets in general). In particular, we prove Hilbert's Basis Theorem. Finally, we investigate bijective correspondences between sets of polynomials and points.
\subsection{Introducing plane curves} 
\begin{defn}
Let $k$ be a field and $f\colon k^2 \rightarrow k$ be a polynomial. A \emph{plane curve} $X$ is defined as the set of all points $\{(x,y)\} \in k^2$ that satisfy $f(x,y)=0$ (such a function is said to ``vanish" on $X$).
\end{defn}
\begin{center} 
\begin{figure}[h!]
  \caption{The plane curve defined by $k=\bbR$ and $f(x,y)=y^2-x^2(x+1)$.}
  \includegraphics[scale=0.1]{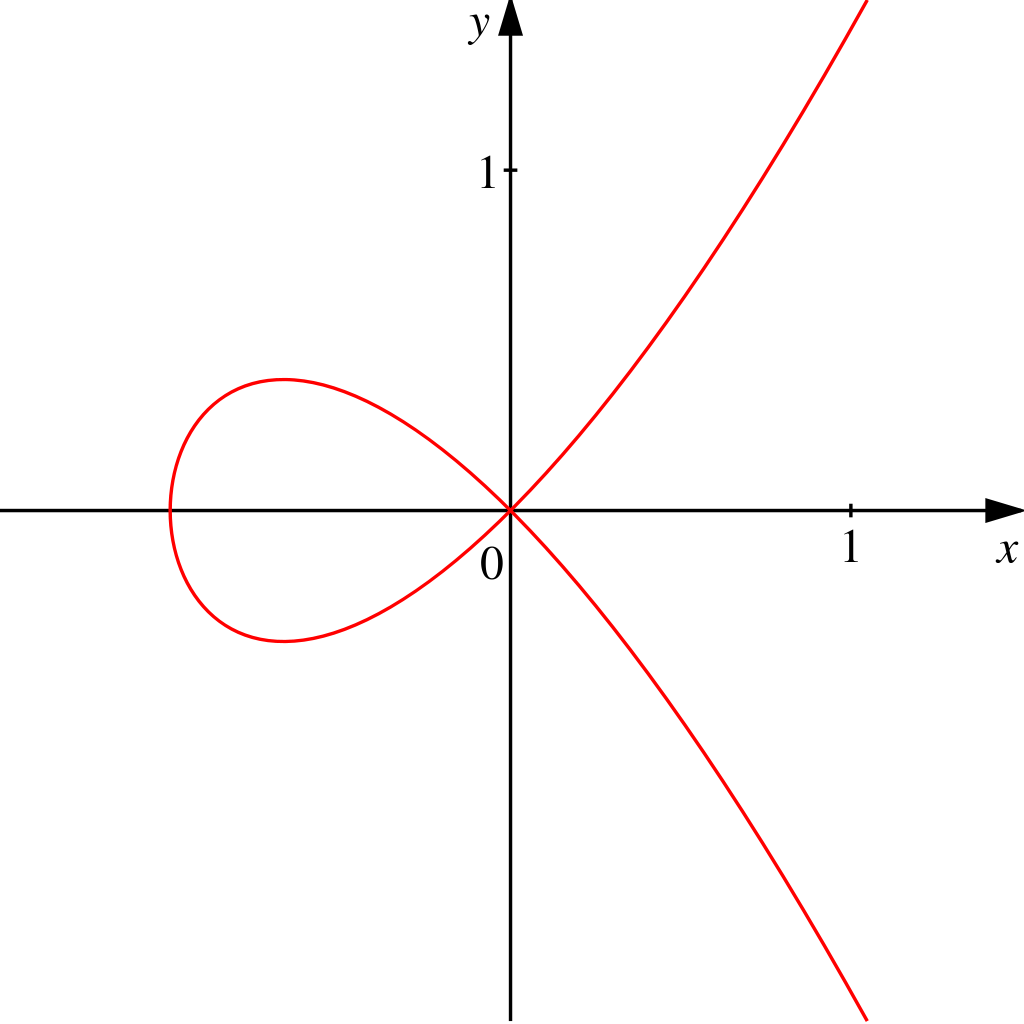}
  \label{fig:planecurve}
\end{figure}
\end{center}
We are usually familiar with finding the zeros of a polynomial $f(x)$ by setting $f(x)=0$, and we know that there are finite number of such $x$ (over a field). Now, we have $f(x,y)$, and there can be an infinite number of $(x,y)$ that satisfy $f(x,y)=0$ (as seen in Figure $\ref{fig:planecurve}$). Algebraic geometry is concerned with studying the relationship between these polynomial functions and sets of points. 
\subsection{The polynomial ring}
We introduce machinery for understanding how to work with these polynomial equations.

\begin{defn}
Let $k$ be a field, $n$ be a natural number, and $A$ be $\bbN^n$. The \emph{polynomial ring} in $n$ variables $k[X_1,\dots,X_n]$ is the set of all expressions of the form $\{\sum_{\alpha}^A a_{\alpha}X^{\alpha}\}$ where $a_\alpha \in k$ for all $\alpha \in A$ and $X^{\alpha}=X_1^{\alpha_1}\dots X_n^{\alpha_n}$ where $\alpha_i$ is the $i$th component in $\alpha$. A \emph{polynomial} in $n$ variables is an element of the polynomial ring. Polynomials are finite sums so we impose the additional restriction that $a_\alpha=0$ for all but finitely many $\alpha \in A$. Let $\{\sum_{\alpha}^A a_{\alpha}X^{\alpha}\}$ and $\{\sum_{\alpha}^A b_{\alpha}X^{\alpha}\}$ be two polynomials. Addition of two polynomials is defined as $\{\sum_{\alpha}^A (a_{\alpha}+b_\alpha)X^{\alpha}\}$. Multiplication is defined as $\{\sum_{\alpha}^A \Pi_\beta^A (a_{\alpha}*b_\alpha)X^{\alpha*\beta}\}$ where $X^{\alpha*\beta}=X_1^{\alpha_1*\beta_1}\dots X_n^{\alpha_n*\beta_n}$. One can prove that the polynomial ring is a ring by Definition \ref{ring}. 
\end{defn}
Note that all concepts we defined about rings and ideals in general apply now to polynomial rings. For example, the polynomial rings has ideals, ring homomorphisms, generating sets, and prime ideals.  
\subsection{The relationship between points and polynomials}
In the beginning of this section, we defined a plane curve. In the last subsection, we generalized the algebraic aspect of plane curves by defining polynomials in $n$ variables. Now, we generalize the ``set of points" aspect of plane curves. 
\begin{defn}
Let $V \colon k[X_1,\dots,X_n] \to k^n$ be a function where
$V(S)=\{x \in k^n \mid f(x)=0 \; \forall \; f \in S\}$. $V(S)$ is known as the \emph{variety} of $S$. $X \subset k^n$ is an \emph{algebraic set} if $X=V(S)$ for some $S \subset k[X_1,\dots,X_n]$. 
\end{defn}
\begin{rem}
Note that we are evaluating $f\in k[X_1,\hdots,X_n]$ on a point $x\in
k^n$. This is done using the field
operations of $k$.
\end{rem}
Under this definition, a plane curve is just an algebraic set created from only one ($f$) two-dimensional (variables $x,y$) polynomial.

We introduce the following proposition that follows from the definition:
\begin{prop}
$\emptyset=V(\{k[X_1,\dots,X_n]\})$ and $k^n=V(\{0\})$ are algebraic sets. 
\end{prop}
When one takes the variety, one maps polynomials to points. Now, we consider a map which takes points to polynomials:
\begin{defn}
Let $X$ be a subset of $k^n$ (not necessarily an algebraic set). The \emph{ideal of} $X$, denoted by $I(X)$, is the set of polynomials in $n$ variables that vanish on $X$. 
\end{defn}
\begin{prop}
$I(X)$ is an ideal.
\begin{proof}
Follows from field operations. 
\end{proof}
\end{prop}
It is clear that the variety of a set of polynomials and the ideal of set of points are somehow fundamentally linked. In which way are these two functions ``inverses" of each other?
\begin{prop}\label{vandi}
The correspondences $V$ and $I$ satisfy the following properties for any field $k$:
\begin{enumerate}
    \item Let $X,Y$ be subsets of $k^n$. If $X \subset Y$, $I(Y) \subset I(X).$
    \item Let $S,T$ be subsets of $k[X_1,\dots,X_n]$. If $S \subset T$, $V(T) \subset V(S)$.
    \item Let $X$ be a subset of $k^n$. Then, $X \subset V(I(X))$. $X=V(I(X))$ if and only if $X$ is an algebraic set.
    \item Let $S$ be a subset of $k[X_1,\dots,X_n]$. Then, $S \subset I(V(S))$.
\end{enumerate}
\begin{proof}
$(1)$ and $(2)$ follow from the definition. 

For $(3)$, note that $I(X)$ is defined as the set of all polynomials that vanish on $X$, so these polynomials will vanish on $x \in X$. For the second part of $(3)$, we know that $X=V(I(X))$; then, $X=V(S)$ for $S=I(X)$, meaning that $X$ is an algebraic set. Next, let us prove that $X$ is an algebraic set implies that $X=V(I(X))$. $X=V(S)$ for some $S$. Suppose, for sake of contradiction, that $I(X)$ does not contain $S$. Then, there exists a polynomial $f \in S, f \not\in I(X)$. We know that $f(x)=0$ for all $x \in X$ by definition, so then $f \in I(X)$---a contradiction. Because $S \subset I(X)$, $V(I(X)) \subset V(S)=X$ so $V(I(X))=X$ (we already know $X \subset V(I(X))$ because of the first part of $(3)$). 

For $(4)$, let $S$ be such a set and let $f \in S$. By definition, we know that $f$ vanishes at the points at which $S$ vanishes. 
\end{proof}
\end{prop}
We figured out under what conditions $X=V(I(X))$ for $(3)$ of the proposition. But we did not figure out under what conditions $S=I(V(S))$. We discuss it later; it is the motivation behind Hilbert's Nullstellensatz. 

Finally, we introduce some operations on elements in the polynomial ring, and what the resulting algebraic set. We hope that these examples provide some intuition on ``taking the variety". 
\begin{prop}
	Let $f_1,f_2 \in k[X_1,\dots,X_n]$. Then,
	\begin{enumerate}
		\item $V(f_1*f_2)=V(f_1)\cup V(f_2)$ and
		\item $V((f_1,f_2))=V(f_1) \cap V(f_2)$.
	\end{enumerate}
	\begin{proof}
		Both follow from the definition.
	\end{proof}
\end{prop}
\subsection{Other properties of algebraic sets}
Before continuing with relationships between points and polynomials, let us come back briefly to algebraic sets. First, consider these two statements:
\begin{prop}
Algebraic sets have the following properties:
\begin{enumerate}
    \item The union of finitely many algebraic sets is an algebraic set.
    \item The intersection of any collection of algebraic sets is an algebraic set.
\end{enumerate}
\begin{proof}
\begin{enumerate}
    \item Let $X=V(S)$ and $Y=V(T)$ be algebraic sets for some sets of polynomials $S$ and $T$. Indeed, $X \cup Y=V(ST)$, where $ST$ denotes the set of all products of an element of $S$ by an element of $T$. Let $x$ be an element in $X \cup Y$. Then, $x \in X$ or $x \in Y$. Without loss of generality, suppose $x \in X$. Then, $f(x)=0$ for all $f \in S$ so every function in $ST$ vanishes on $x$. Now, let $x \in V(ST)$. Without loss of generality, suppose that $x \not\in X$. Then, there is an $f \in S$ such that $f(s) \neq 0$. All polynomials $fg \in ST$, where $g \in T$, vanish on $x$, so then $g(t)=0$ for all $g \in T$ and $x$ vanishes on $Y \subset X \cup Y$.

    By induction, the union of finitely many algebraic sets is an algebraic set as well. 
	\item
    Let $\{X_\alpha=V(S_\alpha)\}$ be a collection of algebraic sets indexed by $\alpha$. Then, $\bigcap_\alpha X_\alpha=V(\bigcup_\alpha X_\alpha)$. Let $x \in \bigcap_\alpha X_\alpha$. Then, $x \in X_\alpha$ for all $\alpha$ so $x \in V(S_\alpha)$ for all $\alpha$. Now, let $x \in V(\bigcup_\alpha X_\alpha)$. By definition, $x \in V(X_\alpha)$ for every $\alpha$ so $x \in X_\alpha$ for every $\alpha$. Then, $x \in \bigcap_\alpha X_\alpha$.
\end{enumerate}
\end{proof} 
\end{prop}
These two statements work for a polynomial in any finite amount of variables, but let us turn our attention back to plane curves. The union of finitely many plane curves is a plane curve, and the intersection of any set of plane curves is a plane curve. 

We now return to the Noetherian property. Is the ring of polynomials in $n$ variables Noetherian? First, we prove the following statement:
\begin{thm}
Let $R$ be a ring. If $R$ is Noetherian, $R[x]$ is Noetherian. 
\begin{proof}\label{noetherian}
Since we work with $R[x]$, any polynomial mentioned in this proof is a polynomial in one dimension.

First, I introduce some definitions: the degree of a polynomial is the smallest possible natural number $i$ such that $a_\alpha=0$ for all $\alpha>i$. In addition, the leading coefficient of a polynomial of degree $i$ is $a_i$. 

Let $I$ be an ideal of $R[x]$. Define $J$ to be the set of leading coefficients of all polynomials in $I$. $J$ is an ideal in $R$:
\begin{enumerate}[label={(\alph*)}]
    \item Closure with respect to addition: let $j,k$ be elements in $J$. They correspond to two polynomials $f,g \in I$. Since $I$ is closed under addition, the polynomial $f+g \in I$. The polynomial $f+g$ has leading coefficient $j+k$ so $j+k \in J$. 
    \item Existence of additive inverse: $0 \in I$ so $0 \in J$. 
    \item If $i \in J_i$ and $r \in R$ then $i*r \in I$: let $i \in J, r \in R$. The element $i$ corresponds to a polynomial $f \in I$. $I$ is an ideal so the polynomial $f*r \in I$. The polynomial $f*r$ has leading coefficient $i*r$ so $i*r \in J$.
\end{enumerate}
$J$ is finitely generated because $R$ is Noetherian. Let $A=\{a_\mu\}$ be the finite set of leading coefficients that generate $J$. For each $a_\mu \in A$, pick $g_\mu \in I$ with leading coefficient $a_\mu$. Let $G=\{g_\mu\}$. By the Well-Ordering Principle, $G$ has a polynomial with largest degree: let $k$ be this degree. 

Then, for every $i \leq k$, define $J^i$ to be the set of leading coefficients of all polynomials in $I$ with degree smaller than or equal to $i$. $J^i$ is an ideal in $R$ (by a similar argument as given above). $J^i$ is therefore also finitely generated. Let $A^i=\{a_{\mu^i}^i\}$ be the finite set of leading coefficients that generate $J^i$. For each $a_\mu^i \in A$, pick $g_\mu^i \in I$ with leading coefficient $a_\mu^i$. Let $G^i=\{g_\mu^i\}$.

Now, let $I_0$ be the ideal generated by all $g^i$ in all $G^i$ and all $g$ in $G$. There are finite number of $g^i$ and $g$ so $I_0$ is finitely generated. I claim now that $I=I_0$. 

Let a polynomial $f \in I_0$. $f$ is generated by some $g$'s and $g^i$'s, which are all elements of $I$, so $f \in I$. 

Now, for sake of contradiction, suppose that $I \not\subset I_0$. Find a polynomial $f$ in $I$ that is not in $I_0$ with smallest degree. We introduce two cases:
\begin{enumerate}[label={Case \arabic*:}, wide=15pt, leftmargin=*]
	\item The degree $\ell$ of $f$ is bigger than $k$. Let $a$ be the leading coefficient of $f$. Then, $a \in J$. This means that $a=\sum_{\mu}^A a_\mu*r_\mu$. Let $f_0=\sum_{\mu}^A r_\mu X^{\text{deg}(f)-\text{deg}(g_\mu)}g_\mu$. The polynomial $f_0$ has the same leading coefficient as $f$ and the same degree as $f$, so $\text{deg}(f-f_0)<\text{deg}(f)$. Also, $f-f_0$ is not in $I_0$ because otherwise by the definition of an ideal, $f$ is in $I_0$. Then, we have found an element with smaller degree in $I$ that is not in $I_0$---a contradiction. 
\item The degree $\ell$ of $f$ is smaller than $k$. Let $a$ be the leading coefficient of $f$. Then, $a \in J_\ell$. This means that $a=\sum_{\mu^i}^{A_\ell} a_{\mu^i}*r_{\mu^i}$. Let $f_0=\sum_{\mu^i}^{A_\ell} r_{\mu^i} X^{\text{deg}(f)-\text{deg}(g_{\mu^i})}g_{\mu^i}$. We get a similar contradiction as the one in Case $1$.
\end{enumerate}
Note that we must separate out the proof into two cases in this way because we do not want the degree of $X$ in $f_0$ to be negative.
\end{proof}
\end{thm}
\begin{cor}\label{hilbert}
$k[X_1,\dots,X_n]$ is Noetherian.
\begin{proof}
$k$ is Noetherian, so $k[X_1]$ is Noetherian. Now, consider $k[X_1,X_2]$. Define $k[X_1][X_2]$ to be the ring of one variable $X_2$ with coefficients in $k[X_1]$. $k[X_1]$ is Noetherian so $k[X_1][X_2]$ is Noetherian. The ring $k[X_1][X_2]$ is equivalent (isomorphic) to the ring $k[X_1,X_2]$. By induction, $k[X_1,\dots,X_n]$ is Noetherian.
\end{proof}
\end{cor}
Theorem \ref{noetherian} is known as Hilbert's Basis Theorem. What can we use it for?
\begin{prop}
Let $X$ be an algebraic set. $X$ is the intersection of the varieties of finitely many polynomials. 
\begin{proof}
$X=V(I(X))$ by Proposition \ref{vandi}. $I(X)$ is finitely generated by Corollary \ref{hilbert}. So, we have $X=V((f_1,\dots,f_n))$. $V((f_1,\dots,f_n))=V(f_1) \cap \dots \cap V(f_n)$ follows from the definition of the variety of a set of polynomials. 
\end{proof}
\end{prop}
This fact is not only interesting when considering algebraic sets (and plane curves), but also useful. In particular, Hilbert used it to prove an important result in invariant theory. 

We use Hilbert's Basis Theorem again. First, we introduce a new definition:
\begin{defn}
An algebraic set $X \subset k^n$ is irreducible if there do not exist $X_1,X_2$ such that $X_1 \cup X_2=X$ with $X_1,X_2 \subsetneq X$.
\end{defn}
We then prove an intermediate result:
\begin{prop}
Every descending chain of algebraic sets $X_1 \supset X_2 \supset \dots$ eventually terminates.
	\begin{proof}
		Suppose that the chain does not terminate. Then, $I(X_1) \subset I(X_2) \subset \dots$ does not terminate---a contradiction because $k[X_1,\dots,X_n]$ is Noetherian.
	\end{proof}
\end{prop}
Finally, we reach the proposition:
\begin{prop}
	Let $X$ be an algebraic set. $X$ can be written as a finite union of irreducible algebraic sets $X_1,\dots,X_n$ where $X_i \not\subset X_j$ for all $i \neq j$ (let this be property $*$). 
	\begin{proof}
		Let $\Sigma$ be the set of all algebraic sets in $k^n$ that do not have property $*$. If $\Sigma$ is empty, our proof is done. Suppose that $\Sigma$ is not empty. Then, there exists an algebraic set $X$ that does not contain any other algebraic set in $\Sigma$ (otherwise we have a descending chain of algebraic sets that does not terminate). $X$ is not irreducible by our given. Then, let $X_1,X_2$ be algebraic subsets such that $X_1,X_2 \subsetneq X$ and $X=X_1 \cup X_2$. One of $X_1,X_2$ cannot have property $*$; otherwise $X$ has property $*$ (by unioning all finite irreducible polynomials making up $X_1$ and $X_2$) so then $X \not\in \Sigma$. Thus, at least one of $X_1,X_2$ is in $\Sigma$, so there does exist an algebraic set in $\Sigma$ that is contained in $X$---a contradiction.
	\end{proof}
\end{prop}
Thus, an algebraic set can be written as the intersection of finitely many varieties of polynomials. In addition, an algebraic set can be written as the union of finitely many irreducible algebraic sets. 

In Section $2$, we did not know of any rings that were not Noetherian. We now give the following example:
\begin{prop}
Let $k$ be a field. $k[X_1,X_2,\dots]$ is not Noetherian.
\begin{proof}
The chain of ideals $(X_1) \subset (X_1, X_2) \subset \dots$ never terminates. 
\end{proof}
\end{prop}
\subsection{Correspondences between points and polynomials}
We introduce a bijective correspondonce from irreducible algebraic sets to prime ideals:
\begin{prop}
Let $X \subset k^n$ be an algebraic set. $X$ is irreducible $\iff$ $I(X)$ is prime.
\begin{proof}
Let us prove that if $X$ is irreducible, $I(X)$ is prime. Suppose that $I(X)$ is not prime. Then, there exist polynomials $f_1,f_2 \not\in I(X)$ such that $f_1f_2 \in I(X)$. Let $I_1$ be the ideal generated by $f$ and all generators of $I$. $I(X) \subsetneq I_1$ so $V(I_1) \subsetneq V(I(X))=X$ because $X$ is an algebraic set. The same logic applies for $f_2$. So we have that $V(I_1),V(I_2) \subsetneq X$. Let $x \in X$. We know that $f_1f_2(x)=0$, so either $f_1(x)=0$ or $f_2(x)=0$, so $x \in V(I_1)$ or $x \in V(I_2)$. Then, $X \subset V(I_1) \cup V(I_2)$ so $X$ is reducible. 

Let us prove now that if $I(X)$ is prime, $X$ is irreducible. Suppose that $X$ is not irreducible. Let $X_1,X_2  \subsetneq X$ such that $X=X_1 \cup X_2$. Because $X_1 \subsetneq X$, $I(X) \subsetneq I(X_1)$, there exists $f_1$ such that $f_1 \in I(X_1)$ but $f_1 \not\in I(X)$ (and a similar $f_2$ for $I(X_2)$). $f_1f_2$ vanishes at all points in $X$ so $f_1f_2 \in I(X)$, so $I(X)$ is not prime. 
\end{proof}
\end{prop}
Now we return to Proposition \ref{vandi} $(4)$ to see if we can find conditions such that $S=I(V(S))$. If we succeed, we can define a bijective correspondence between these sets $S$ and algebraic sets $V(S)$.  Why are we not guaranteed equality?

Instead of considering all sets of polynomials $S$, we consider ideals of the polynomial ring $J$ only in this problem. Recall that if a set of polynomials $S$ and an ideal $J$ both vanish on $X$ (i.e. $V(S)=V(J)$), $S \subset J$. As we are trying to achieve $S=I(V(S))$, using $S \neq J$ would mean that $S \subset I(V(S))$---not what we want.  

First, consider the polynomial $x^2+1 \in k[x]$, where $k=\bbR$. Let $J$ be the ideal generated by this polynomial. $V(J)=\emptyset$, so $I(V(J))=k^n$. As we can see, the solution just doesn't exist in $\bbR$. If we were to use  $k=\mathbb{C}$, though, $V(J)$ would not be empty (actually $|V(J)|=2$). So, in order for equality to hold, we need $k$ to have some special property. 

Second, consider the polynomial $x^2 \in k[x]$, where $k=\bbR$. Let $J$ be the ideal generated by this polynomial. $V(J)=\{0\}$, and $I(V(I))=(x)$. As we can see, the original ideal did not contain all polynomials that vanish on its vanishing set, so we ended up with inequality. Then, if we want $J=I(V(J))$, for every $f^n$ in $J$, $f$ must be in $J$. 

We define two new definitions based on these two observations:
\begin{defn}
A field $k$ is algebraically closed if every polynomial in $k[x]$ has a root.
\end{defn}
\begin{defn}
Let $J$ be an ideal. The \emph{radical} of $J$, denoted by $\sqrt{J}$, is the set $\{f \mid f^n \in J \text{ for some } n>0\}$. $\sqrt{J}$ is an ideal. 
\end{defn}
Hilbert's Nullstellensatz proves the bijective correspondence between radical ideals and algebraic sets: 
\begin{thm}
Let $k$ be an algebraically closed field and $J$ be an ideal with coefficients in $k$. Then, $\sqrt{J}=I(V(J))$. 
\end{thm}
Proof of Hilbert's Nullstellensatz requires more definitions; we do not prove it here. We hope though that the reader has noted how useful such statements are to understanding the polynomial ring and these sets of points. 

\section*{Acknowledgements}
I am very grateful to my mentor Michael Neaton and to Dr. May for their useful comments and suggestions. All errors are mine.

\end{document}